\documentclass{article}
\usepackage{epsfig,amssymb}
\begin{document}
\newtheorem{Th}{\sc Theorem}
\newtheorem{Lm}{\sc Lemma}
\newtheorem{Df}{Definition}
\newcommand{\eqdef}{\stackrel{def}{=}}
\newcommand{\rlm}{\rho_{(\lambda,\mu)}}
\newcommand{\wh}{\widehat}
\newcommand{\tomeg}{\tilde\omega}
\newcommand{\tkap}{\tilde\kappa}
\newcommand{\tzet}{\tilde\zeta}
\newcommand{\Ob}{{\cal O}}
\newcommand{\F}{{\cal F}}
\newcommand{\Hb}{{\cal H}}
\newcommand{\tP}{\tilde P}
\newcommand{\der}{\mathop{\rm d}}
\newcommand{\mod}{\mathop{\rm mod}}
\newcommand{\ddt}{\frac{d}{dt}}
\newcommand{\HSf}{\frac{\partial{\cal H}}{\partial S_1}}
\newcommand{\HSs}{\frac{\partial{\cal H}}{\partial S_2}}
\newcommand{\Sff}{\frac{\partial S_1}{\partial f}}
\newcommand{\Ssf}{\frac{\partial S_2}{\partial f}}
\newcommand{\SfF}{\frac{\partial S_1}{\partial F}}
\newcommand{\SsF}{\frac{\partial S_2}{\partial F}}
\newcommand{\pF}{\partial_F}
\newcommand{\pH}{\partial_H}
\newcommand{\Tpf}{\frac{\partial T}{\partial\phi_1}}
\newcommand{\pf}{\frac{\partial}{\partial\phi_1}}
\newcommand{\ps}{\frac{\partial}{\partial\phi_2}}
\newcommand{\pT}{\frac{\partial}{\partial T}}
\newcommand{\ppf}{\frac{\partial}{\partial f}}
\newcommand{\FP}{{\cal F}_{P_c}}
\newcommand{\FQ}{{\cal F}_{Q^3_c}}
\newcommand{\FU}{{\cal F}_{U(Q^3_c)}}
\newcommand{\sgr}{\mathop{\rm sgrad}}
\newcommand{\la}{\lambda}
\newcommand{\al}{\alpha}
\newcommand{\be}{\beta}
\newcommand{\ga}{\gamma}
\newcommand{\de}{\delta}
\newcommand{\ka}{\kappa}
\newcommand{\io}{\iota}
\newcommand{\om}{\omega}
\newcommand{\p}{\partial}
\newcommand{\fr}{\frac}
\newcommand{\ep}{\varepsilon}
\newcommand{\ddf}{\frac{d}{d\phi}|_{\phi=0}}
\newcommand{\ct}{\mathop{\small\rm const}}

\title{\Large  Jacobi vector fields of integrable geodesic flows
 }
\author{Vladimir  Matveev\footnote{vmatveev@physik.uni-bremen.de}, 
Peter  Topalov\footnote{topalov@banmatpc.math.acad.bg} }
\date{}
\maketitle

\begin{abstract}

We show that an invariant surface allows to construct the Jacobi vector
field along a geodesic  and construct
the formula for the normal
component  of the Jacobi field. If a  geodesic
 is the transversal
intersection  of two invariant surfaces (such situation we have,
 for  example, if the geodesic  is  hyperbolic),
  then we can construct
a  fundamental solution of the the Jacobi-Hill equation $\ddot u=-K(u)u$.
This is  done for quadratically integrable geodesic flows.
\end{abstract}

\centerline{\sc \S1. Introduction. }

\vspace{.2cm}
{\bf 1.1. Definitions.}   
Suppose  $G=(g_{ij})$  is a Riemanian  metric on a  surface  $P^2$,
a curve  \newline $\gamma:[a,b]\to P^2$ is a geodesic.
We will assume that the  parameter $t\in [a,b]$ of the geodesic $\gamma$
 is natural or natural,  multiplied by a constant.

\begin{Df}
{\em Geodesic variation  } of a geodesic  $\gamma$
is called the  smooth
mapping $\Gamma:[-\ep,\ep]\times [a,b]\to P^2$ such that
\begin{itemize}
\item[{\bf 1)}]  for any fixed  $s_0\in [-\epsilon,\epsilon]$
the curve $\Gamma(s_0,t):[a,b]\to P^2$ (as the curve of parameter
  $t\in[a,b]$) is a geodesic,
\item[{\bf 2)}] for any $t\in [a,b]$  \ $\Gamma(0,t)=\gamma(t)$.
\end{itemize}
\end{Df}

\begin{Df}
{\em  Jacobi vector field }
along the geodesic  $\gamma$ is
the vector field
$J=\frac{\partial\Gamma}{\partial s}{|_{_{s=0}}}$, where
$\Gamma$ is a geodesic
variation of the geodesic  $\gamma$.
\end{Df}

By definition,   Jacobi vector field is a smooth vector field along the
 geodesic.

\begin{Df}
          Jacobi vector field  $J$ is called
          {\em normal} if it is orthogonal
to the geodesic  at every point of the geodesic.
\end{Df}

It is known that the projection of a Jacobi field $J$ to the vector field of
 normals to the geodesic is a normal Jacobi vector field.

 The length of a normal Jacobi vector field  $J$ satisfies
 the {\em
 Jacobi-Hills equation
for   the normal component} \
$\ddot x+K(\gamma(t))x=0$,
where  $K$ is the Gauss curvature
 and  $t$ is the natural parameter (see, for example,
 \cite{DNF} or \cite{Bes}).

Consider real numbers $a\ne b$. Denote by  $A$ the point
 $\gamma(a)\in P^2$, denote by  $B$ the point
  $\gamma(b)\in P^2$.

  \begin{Df}
  The points $A$ and $B$ are called {\em conjuate}
   along the  geodesic $\gamma$
  if  there exists a non-zero Jacobi vector field $J$ along
  the geodesic
  $\gamma$ such that  $J(a)=J(b)=0$. (Figure 1)
  \end{Df}

\begin{figure}[h]
{\psfig{figure=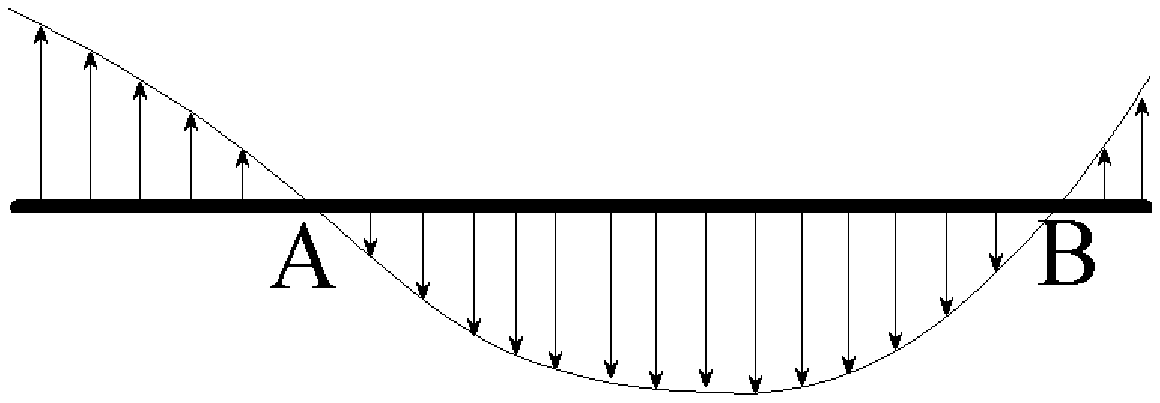}}
\caption[]{\footnotesize }
\end{figure}

%

The point  $A$ can coincide with the point
$B$. It happens if the geodesic  $\gamma$ is closed
or self-intersecting. In the first case the point $A$  is called {\em
self-conjugate } along the
 geodesic  $\gamma$.

\vspace{3mm}
{\bf 1.2. Jacobi vector fields as the projection of  invariant
vector fields  from the co-tangent space.}
The metric allows to identify canonically the tangent
 and co-tangent bundles of the
surface
$P^2$. Therefore we have a scalar product and a norm on every
co-tangent plane. For example,
suppose   $G$
in coordinates $(x,y)$
reads  $\lambda(x,y)(dx^2+dy^2)$. Then the scalar product on
  $T^*P^2$ is given by the formula
$<(p_x,p_y),(\hat p_x,\hat p_y)>=\frac{p_x\hat p_x+p_y\hat p_y}{\lambda(x,y)}$.

\begin{Df}
{\em Geodesic flow of the metric  $G$} is called the Hamiltonian system
on $T^*P^2$ with
the Hamiltonian $H\eqdef\frac{1}{2}|p|^2$, where  $p$ is  momentum and
$|.|$ is the norm.
\end{Df}

In particular, the Hamiltonian $H$ of the geodesic flow of the metric
 $\lambda(x,y)(dx^2+dy^2)$   is given by the formula
  $H(x,y,p_x,p_y)=\frac{p_x^2+p_y^2}{2\lambda(x,y)}$.

It is known that  the  trajectories of the geodesic
 flow projects onto  the geodesics.

 \begin{Df}
An imbeded into  $T^*P^2$ surface $I^2$
is called {\em invariant}
 if the vector field of the geodesic flow is tangent to $I^2$.
\end{Df}

 \begin{Df}
Let $\gamma:[a,b]\to T^*P^2$  be the solution of the geodesic flow.
The vector field $V$ along the curve $\gamma$ is called {\em invariant}
 if
it is invariant with respect to the  the geodesic flow.
 \end{Df}

In other words, consider the one parametric family
of mappings $S_{\tau}:T^*P^2\to T^*P^2$. The mapping
$S_{\tau}$ moves a point
$x\in \ga\subset  T^*P^2$  along the trajectory $\ga$ of the geodesic flow
during  the time $\tau$.
The vector field $V$ is called {\em  invariant}
 if for any $\tau$ the differential
$dS_{{\tau}_{|_x}}$ takes  the vector field $V$  to itself.

 In \S2 we will show that a
geodesic variation allows to canonically construct
an invariant surface;
a Jacobi vector field  allows to canonically construct
an invariant vector field;
the   projection of the
invariant vector field is the Jacobi vector field and
the composition of the projection and an
 imbedding of the invarant surface is the  geodesic variation.

\vspace{2mm}

{\bf 1.3.  Jacobi vector fields of integrable geodesic flows.
}
A geodesic flow is
called integrable if it
is integrable as a  Hamiltonian system.
That is there exists a smooth function
 $F:T^*P^2\to R$ such that
\begin{itemize}{
\small  \item[{\bf 1)}] $F$
is constant on the trajectories of the geodesic flow
 \item[{\bf 2)}] the differentials  $dH$ and $dF$ are linear independent
 almost everywhere.
                }
\end{itemize}

The function  $F$ is called  {\em an  integral.}
Note that the geodesic flow preserves  the vector field
$\mbox{sgrad}(F)$.
Since that, the projection of the vector field
 $\mbox{sgrad}(F)$ is a Jacobi vector field.
 Using this, we can  construct a number of pairs of conjugate points.
Let  $L^2$ be a Louville torus of an  integrable geodesic flow.
 Restrict the natural projection to the torus $L^2$.
 A connected component of the set of
 critical points of $\pi_{|_{L^2}}$ is called
 a caustic. It is known  ( see, for example, \cite{BP})
 that a caustic
 is a smooth simple curve and  can not intersect other caustic.

{\bf Remark}  {\small   Sometimes caustics are called
                         the projection  of the set of
                          critical points of the mapping
                          $\pi_{|_{L^2}}$.
                          According to our definition,
                           caustics are curves in the
                            fase space $T^*P^2$.
                                                }

Suppose  the intersection of the trajectory
 $\hat\gamma$ with the caustics includes
  points
 $\hat A\in L^2\subset  T^*P^2$ and $\hat B\in L^2\subset  T^*P^2$.
We prove that  the projections
$\pi(\hat A)$ and
$\pi(\hat B)$ are conjugate along the
geodesic  $\gamma=\pi(\hat\gamma)$. Indeed,
consider the restriction of
 the vector field
$\mbox{sgrad}(F)$ to the Liouville
torus  $L^2$. Consider
the projection
 $d\pi(\mbox{sgrad}(F))$.
 Since the rank of
the projection
 $d\pi(\mbox{sgrad}(F))$ is less than 2 in the points
             $\pi(\hat A)$ \ and  \ $\pi(\hat B)$, we see that
the projections of the vector fields
 $\mbox{sgrad}(H)$ and
$\mbox{sgrad}(F)$ are parallel in the points.
Therefore the normal component of the vector field $\pi \mbox{sgrad}(F)$
equals zero in the points $\pi(\hat A)$ and $\pi(\hat B)$. Thus points
$\pi(\hat A)$ and $\pi(\hat B)$ are conjugte.

 \vspace{2mm}
{\bf 1.4. Jacobi vector fields for a  hyperbolic geodesic.}
\begin{Df} A closed geodesic $\gamma$ is
called {em hyperbolic} if
the corresponding trajectory $\hat\gamma$ of
the geodesic flow is hyperbolic.
 \end{Df}

That is if we restrict the geodesic flow to the
isoenergy surface (we denote isoenergy surface by $Q^3$)
$\{H=\mbox{\small const}\}$ then the multipliers
 of
the trajectory   not lie
on the
unit circle.

 It is known that
 in a regular neighborhood $U(\hat \gamma)\subset Q^3$  of a
hyperbolic trajectory there exists
a pair of invariant  two-dimensional surfaces
($L_+$ and $L_-$). The intersection of the
surfaces coincides with $\hat \gamma$
(See, for example, \cite{Ano}).

An invariant surface allows to construct
a solution of the Jacobi-Hill equation.
In  \S3 we will show that the
solutions that correspond to the surfaces $L_+$ and
$L_-$  are
not  proportional.
Therefore we have  a fundamental solution
of the Jacobi-Hill equation  $\ddot x+K(\gamma(t))x=0$.

Assume that the geodesic flow is integrable, and the integral is a Bott function.
That is the restriction of the integral to the isoenergy
surface satisfies the following properties (we denote the restriction by $F$):
{\small \begin{itemize}
\item[{\bf 1.}]   The critical manifolds of  $F$ are compact sets.
\item[{\bf 2.}]    If $D$ is an arbitrary
 2-disk that is transversal with  critical manifolds,
 then the restriction $F_{|_{D}}$
 is a Morse function.
        \end{itemize} }

Suppose  a connected component of the set of
 critical points is homeomorphic to
the circle. Consider   a transversal disk.
The
 dimension of the transversal disk
equals 2. The restriction of the function $F$
to the transversal disk
has Morse singularity of index 0,2, or {\bf  1}.
In the last  case the connected component
of the critical  set is called {\it a saddle circle}.
Hyperbolic trajectories of  a  Bott integrable geodesic flow
are saddle circles.

\vspace{.2cm}
{\bf Remark } { \small Let a  trajectory of a
Bott integrable
Hamiltonian system  be
 a hyperbolic trajectory.
Then it is a saddle circle.
The reverse statement is not true.
There exist saddle circles that are
not hyperbolic trajectories. }

\vspace{2mm}
Consider the Liouville fiber that contains a saddle circle.
In a  neighborhood of a point of the saddle circle the Liouville fiber is
homeomorphic to a pair of intersecting  surfaces.

Recall that   a saddle circle is called
{\it nonorientable} if the intersection
of the Liouille
fiber with  a regular neigborhood of the saddle circle
is homeomorphic to the self-intersecting
 M\"obius band.
A saddle circle is called {\it orientable}
the intersection of  the Liouville fiber    with a saddle neighborhood
is homeomorphic to two intersecting annuli.

Consider a closed  geodesic  $\ga$.
Let a point $B\in \ga$ be self-conjugate
along the geodesic. Consider the set of  points
conjugate to $B$.
Denote by $N(\ga)$ the number of  elements in the set.
It is known that $N(\ga)$ does not depend on the
choice of the initial point $B\in\ga$.

If for the point  $C\in \ga_1$ there are no conjugate  points
 along $\ga_1$, then
by defintion put $N(\ga_1)=0$.

\begin{Th} Suppose  $P^2$ is an orientable surface,
the geodesic flow of a metric
$G$ is integrable, and     $\hat\gamma$ is a saddle circle.
Then
\begin{itemize}
\item if the circle  $\hat\ga$ is nonorientable, then $N(\pi(\hat{\ga}))$
 is odd number,
\item  if the circle $\hat\ga$ is orientable,
 then $N(\pi(\hat\ga))$
is even nuber. \end{itemize}
\end{Th}

The statement of the result for nonorientable saddle circles
were proved by
A. Wittek.

                         \vspace{.2cm}

                         \vspace{.2cm}
Proof. Let a saddle circle  $\hat\ga$ be  orientable.
Consider the Liouville fiber that contains $\hat\ga$.
 By definition,  a regular
 neighborhood of $\hat\ga$ the Liouville  fiber
is homeomorphic to a pair of intersecting annuli.
Denote by $L_+$  oe of the    annuli.
Consider a
vector field $\hat\ga$ that is invariant and that is tangent
 to
 $L_+$.
  Denote by $J$ the projection $\pi(\hat\ga)$ of the vector field
 $\hat J$.

Using the theorem of the existence and
uniqueness
of the solution of a differential equation, we have
that in a
neighborhood of a zero point a normal Jacobi vector
field behaves as it is shown on
Diagramm  2(a) (the situation of Diagramm 2(b) is forbidden).
That is,  the frame  (velocity vector of the geodesic,
Jacobi vector field)
has different orientation at different sides of the
geodesic $\ga$.

\begin{figure}[h]
  \centerline{\psfig{figure=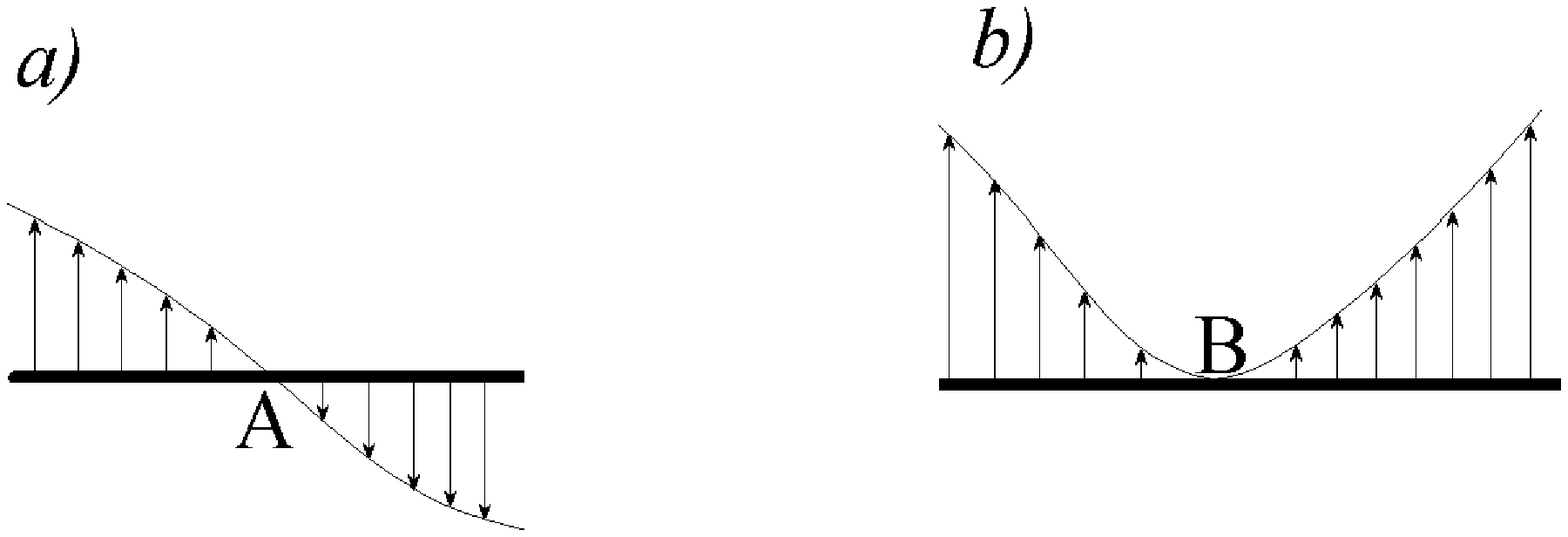}}
\caption[]{\footnotesize }
\end{figure}

%

Since  $\hat J$ has no zero points, we see that the frame
(the velocity vector of  $\hat \ga$,   $\hat J$)
has the same orientation with respect
to the invariant surface $L$ in every point of the geodesic.
Hence the vector field $J$ has the same direction after and before
the circuit along the geodesic. Hence, the number
of zeros of $J$ is an even number.
Proved.

\vspace{2mm}

{\bf 1.6. Jacobi vector fields along hyperbolic geodesics
 of quadratically integrable geodesic flows.}

\begin{Df}A geodesic flow is called  {\em
linear integrable } if it admits an integral $F$
such that in a neighborhood of any point the integral $F$
is given by the formula
 $F(x,y,p_x,p_y)=A(x,y)p_x+B(x,y)p_y$, where  $x,y$ are  coordinates
 on the surface,  $p_x,p_y$
  are the corresponding  momenta,
  and  $A,B$ are smooth functions of two variables.

\end{Df}

\begin{Df}
A geodesic flow is called
 {\it quadratically integrable},
 if it is not linear integrable and if
 it admits an integral $F$ such that in a neighborhood of any
  point the integral $F$ is given by the formula
 $F(x,y,p_x,p_y)=A(x,y)p_x^2+B(x,y)p_xp_y+C(x,y)p_y^2$, where
  $A,B,C$ are smooth
 functions of two variables.
\end{Df}

Quadratically and linear integrable geodesic flows on closed surfaces
are comletely described. In \cite{Koz} V.V. Kozlov proved that there are no
linear and quadratic integrable geodesic flows on the surfaces of genus $g>1$.
  Linear and quadratically integrable geodesic flows on the sphere
were described by V.N. Kolokolzov in   \cite{Kol}.
Quadratically integrable geodesic flows on the torus were described
by I.K. Babenko and N.N. Nekhoroshev in \cite{BN} (see, also, \cite{Mat1}).

In \S5 we describe hyperbolic trajectories
and invariant surfaces of
quadratically inegrable geodesic flows and obtain
 fundamental solutions  of the Jacobi-Hill equation
 $\ddot x+K(\gamma(t))x=0$ along hyperbolic geodesics.

 \vspace{3mm}

 \centerline{\sc \S2. Canonical frame on $T^*T^2$. }
\centerline{\sc Commutative relations for it.}

\vspace{.2cm}

Let $G$ be a Riemannian metric on the oriented surface  $P^2$.
 Consider the tangent bundle  $TP^2$.
  The space of non-zero vectors of  $TP^2$
  is denoted by  $T_0P^2$.

  The aim of this section is to canonically construct
  a frame (the vectors of the
  frame will be denoted by $D_\phi, D_1, D_2,A$) in
  every point of
   $T_0P^2$.  Since the metric allows to identify the
   tangent and co-tangent
   bundles of the surface $P^2$, we will
   have the canonic frame in every point of $T_0^*P^2$.

  {\bf 2.1. The frame. } Let  $x_0$ be an arbitrary point
  of the surface  $P^2$,
    let $v$ be a tangent vector at the point $x_0$.
     The vector $v$ is a point of  $T^*P^2$.

    Denote by   $\rho_{\phi}(v)$ the vector $v_1$ at the point $x_0$
    such that
    $|v_1|=|v|$, the angle between $v$ and $v_1$ is equal to $\phi$, and
    the frame   $(v, v_1)$ is positive.  In other words,
     $\rho_{\phi}$
     rotates  a vector by     the  angle $\phi$.
  Consider the vector
  $D_{\phi}\eqdef\ddf\rho_{\phi}(v)$ at the point $v\in TP^2$.

  Denote by  $D_1$ the vector field $\mbox{sgrad}(H)$.

  By definition, put  $
  D_2\eqdef[D_{\phi},D_1].$

  Consider the so-called "Liouville vector field" $A$ (\cite{Bes}).
  Recall, that       the vector field $A$ is defined in the following way.
  Consider the one-parameter group of self-diffeomorphisms $g^t:TP^2\to TP^2$.\,
    $g^t(v)\eqdef exp(t)\vec v$.  Put by definition
  $A(v)\eqdef\fr{d}{dt}|_{t=0}g^t(v)$.

      Direct calculations show  that for $v\in T_0P^2$ the
      vectors $D_\phi, D_1, D_2, A$ are linear independent.
      Therefore the quadruple  $(D_\phi, D_1, D_2, A)$ is a frame.

\vspace{.2cm}
{\bf 2.2. Commutative relations for the vectors of the frame.}
Denote by  $r(v)$  the function  $\sqrt{G(v,v)}$. The function
 $r(v)$ is a smooth function on
 $T_0P$.

  \begin{Lm}
   $$
\left\{
  \begin{array}{l}
  [D_1,D_2]=r^2KD_{\phi},\cr
  [D_{\phi},D_1]=D_2,\cr
  [D_{\phi},D_2]=-D_1,\cr
  [A,D_1]=D_1,\cr
  [A,D_2]=D_2,\cr
  [A,D_{\phi}]=0,\cr
  \end{array}
\right.
 $$
  where $K(v)\eqdef K(\pi(v))$,  $K;P2\to R$ is the Gaussian
   curvature of the metric
   $G$ on the surface  $P^2$.  \end{Lm}
The first three relations follows from \cite{KN}.
The last  three relations were proved in  \cite{Bes}.

Consider a geodesic trajectory  $\hat\ga(t)$. Consider an invariant
vector field
 $\hat J$ along $\hat\ga(t)$. Let $
\hat J=xD_2+yD_{\phi}+ID_1+aA.
 $

Since the vector field  $\hat J$ is invariant, we have $[D_1,\hat J]=0$.
Using lemma 1, we obtain the following system:
       \begin{equation}\label{1}
 \cases{\dot x=y,\cr
       \dot y+r^2Kx=0,\cr
        \dot a=0,\cr
       \dot I=a.\cr}
   \end{equation}

If the length of the vector  $\pi(D_1)$ equals  1, then from
the first two equations  it follows Jacobi-Hills
 equation  $\ddot x+Kx=0$.
Note, that in this case
$x$                     is normal component
of the vector field $\pi(\hat J)$.
Projection of the vector field  $I$ is the horizontal component
of the vector field  $\hat J$. Hence the projection $J$ of the vector field
 $\hat J$ is equal to  $x\vec n+a\dot\gamma$, where
$\vec n$ is the normal vector (to the geodesic) of length 1.

\vspace{2mm}

We shall prove that for every Jacobi vector field  $J$ there exists
an  invariant vector field  $\hat J$ such that  $\pi(\hat J)=J$.

Consider a Jacobi vector field  $J(t)$ along the geodesic $\ga(t)$.
Denote by
$x(t)$ the normal component of the vector field $J$, denote   by $I(t)$
the horizontal component of the vector field $J$.
By definition, $J(t)=\fr{\p}{\p s}|_{s=0}\Gamma(s,t)$, where
 $\Gamma(s,t)$ is a geodesic variation of the geodesic
 $\dot\ga$. Consider
 the vector field  $\hat
J\eqdef\fr{\p}{\p s}|_{s=0}\dot\Gamma(s,t)$ along the
 trajectory  $\hat\ga$.
Evidently the vector field $\hat J$ is invariant. Indeed,
 $S^\tau(\dot\Gamma(s,t))=\dot\Gamma(s,t+\tau)$.

Moreover, we see that the vector  field
 $\hat J$ equals   $$ \hat
J=x(t)D_2+y(t)D_{\phi}+I(t)D_1+a(t)A, $$ where $y(t)$, $a(t)$ are smooth functions,
 $x$ is the normal component of the Jacobi vector field $J$, and
  $I$ is the horizontal component of
the Jacobi vector field  $J$.

 \vspace{2mm}

 \vspace{3mm}

 \centerline{\sc \S3. Fundamental solution of the
 Jacobi-Hill  equation}
 \centerline{\sc for hyperbolic geodesics.  }

 \vspace{2mm}
 Consider a geodesic trajectory  $\hat\ga(t)$.
 Denote by  $\ga$ the geodesic  $\pi(\hat\ga(t))$.

  Suppose we have an 1-dimension subspace in every tangent
  space at   the points  of the
  geodesic trajectory.
  Suppose  the set of the subspaces $l_{inv}$
  is invariant with respect ot the geodesic
   flow.

Consider the projection of every subspace $l_{inv}$ along the plane
    $<A,D_1>$ to the plane $<D_2, D_\phi>$. Denote the projection by
    $l'_{inv}$.

Consider the mapping  $d\pi:T(TP^2)\to TP^2$. Then
 $d\pi(l'_{inv})$ is 1-dimensional
 subspase in the tangent planes  at
 the points of the geodesic.  Using Lemma~1, we see that
  the set of $(l'_{inv})$  is invariant with respect to
   the geodesic flow.
 Since the vector $D_\phi$ projects to $0$, and since
 the projection      of $D_2$
   is orthogonal to the geodesic  $\ga$, we have that
  $d\pi(l_{inv})$  is orthogonal to the geodesic  $\ga$.

  Consider a direction vector field $Y$
   of the set $l'_{inv}$.
  Let an invariant  vector field $\xi$ be parallel to $Y$.
 Then  $\xi(t)$  equals  $\ka(t)Y(t)$,
 where $\ka(t)$ is
 a smooth  function.
 Suppose   the coordinates of the vector field
     $Y$ in the frame $<D_\phi,D_1,D_2,A>$  are equal to
     $(\al(t),0,\be(t),0)$.
      (Since $Y$ lies in the plane $<D_\phi,D_2>$, we see that
     the second and the fourth coordinates must  be equal to $0$.)

We have  $Y(t)=\al(t)D_2+\be(t)D_{\phi}$.  Then,
 $\xi(t)=\ka(t)\al(t)D_2+\ka(t)\be(t)D_{\phi}$.
 Using  (\ref{1}), we obtain
 $ \fr{d}{dt}{(\ka\al)}=\ka\be.$
 Therefore,
   $$
 \ka(t)=\exp{\left[{\int_{t_0}^t\fr{\beta(s)-\dot\al(s)}{\al(s)}ds}\right]}=
\fr{\al(t_0)}{\al(t)}\exp\left[\int_{t_0}^t\fr{\beta(s)}{\al(s)}ds\right].
  $$
  Thus the vector field
  \begin{equation}\label{2}
 \xi=\fr{\al(t_0)}{\al(t)}\exp\left[\int_{t_0}^t\fr{\beta(s)}
{\al(s)}ds\right]Y \end{equation}
  is invariant with respect to the geodesic flow, and the projection
  of it is a Jacobi vector field.

 In particular, the function
 \begin{equation}\label{3}
     \exp\left[{\int_{t_0}^t\fr{\beta(s)}{\al(s)}}ds\right]
     \end{equation} satisfies the Jacobi-Hill equation
      $\ddot x+r^2Kx=0$.

 Suppose a geodesic  $\hat\ga$
 is hyperbolic. Then we have
 two invariant intersecting surfaces in a regular neighborhood of $\hat\ga$.
 Consider the subspaces  that are tangent to  these surfaces.
 Since the intersecting  surfaces are
 invariant then the tangent subspaces
 are invariant, too. Then the projection of the surfaces to the planes
 $<D_\phi,D_2>$
 along the planes
 $<D_1,A>$ is two families of   invariant 1-dimension subspace.
 Arguing as above, we can construct
 two Jacobi vector fields.  The Jacobi vector field are nonproportional.
 Hence they  define  a fundamental solution of the Jacobi-Hill equation.

  {\bf Remark }  {\small  It is easy to find the coordinates  $\al,\beta$.
  Let the metric $G$ has the form
 $\lambda(dx^2+dy^2)$.   Let a field $Y$ equals
  $k_1\frac{\p}{\p x}+k_2\fr{\p}{\p y}+K_1\fr{\p}{\p
 p_x}+K_2\fr{\p}{\p p_y}$.  The coefficients
  $k_1,k_2,K_1,K_2$ are functions of  $t$.  Then,
 \begin{equation}\label{4}
 \al=\lambda\fr{\dot xk_2-\dot yk_1}{r^2},
  \end{equation}
 \begin{equation}\label{5}
  \beta=\left(\fr{1}{2\la}\fr{\p\la}{\p x}k_2
 -\fr{1}{2\la}\fr{\p\la}{\p y}k_1\right)+
 \fr{\dot xK_2-\dot yK_1}{r^2}.
 \end{equation}
                  }

  {\bf Remark } {\small
Now let the projection  $d\pi(l_{inv})$  of an invariant set of subspaces
$l_{inv}$ is not orthogonal to the geodesic.
Consider the projection of every subspace $l_{inv}$ along the plane
    $<A,D_1>$ to the plane $<D_2, D_\phi>$. Denote the projection by
    $l'_{inv}$.
Let  $Y$ be a direction  vector field of the subspaces  $l_{inv}$, denote by
 $Y'$ the
 projection of the vector field $Y$  along the plane
    $<A,D_1>$ to the plane $<D_2, D_\phi>$.
 Suppose that in the coordinate system   $(x,y,p_x,p_y)$ vector fields
$Y$, $Y'$  are equal to  $(k_1,k_2,K_1,K_2)$,
$(k_1',k_2',K_1',K_2')$, respectively. Then,
 $$\lambda\fr{\dot xk_2-\dot yk_1}{r^2}
  =\lambda\fr{\dot xk_2'-\dot yk_1'}{r^2},
  $$
 $$
  \left(\fr{1}{2\la}\fr{\p\la}{\p x}k_2
 -\fr{1}{2\la}\fr{\p\la}{\p y}k_1\right)+
 \fr{\dot xK_2-\dot yK_1}{r^2}=
  \left(\fr{1}{2\la}\fr{\p\la}{\p x}k_2'
 -\fr{1}{2\la}\fr{\p\la}{\p y}k_1'\right)+
 \fr{\dot xK_2'-\dot yK_1'}{r^2}.$$

 }

 \vspace{.2cm}

 Now let the curve  $\ga(\tau)=(x(\tau),y(\tau))$ is a geodesic.
 We do not
 require the parameter  $\tau$ to be natural or natural,
  multipling  by a constant. Let the cordinates
  $(x,y)$ be isothermal.
  That is the metric $G$ has the form
   $\la(x,y)(dx^2+dy^2)$. Suppose  a geodesic trajectory
    $\hat\ga$   projects into the geodesic
     $\ga$. Let $r$ be the length of the
     velocity vector of the geodesic in the parameter $t$.
     Consider the functions   $k_1,k_2,K_1,K_2$ as the functions of the
     new parameter  $\tau$.

  \begin{Lm} \label{lemma}   Let $\tau(t)$ be the connection between the
  parameter  $\tau$ and the parameter $t$.
  Consider the function
$u(\tau)=\exp\left[\int\limits_{\tau_0}^\tau\fr{\beta(s)}{\al(s)}\,|\fr{d\ga}{ds}| ds \right] , \ \ \ \mbox{where}$ $$
 \al(\tau)=\lambda(\tau)\fr{\fr{dx}{d\tau}k_2-
 \fr{dy}{d\tau}k_1}{r^2|\fr{d\ga}{d\tau}|}, \
 \beta(\tau)=\left(\fr{1}{2\la(\tau)}\fr{\p\la}{\p x}k_2
 -\fr{1}{2\la}\fr{\p\la}{\p y}k_1\right)+
 \fr{\fr{dx}{d\tau}K_2-\fr{dy}{d\tau}K_1}{r^2|\fr{d\gamma}{d\tau}|}.
 $$ Then the function $u(t)\eqdef u(\tau(t))$  is
 a solution of the Jacobi-Hill equation.  \end{Lm}

The proof is by direct calculation.

         \vspace{3mm}

  {\sc\centerline{ \S4. Quadratically integrble geodesic flows }
   \centerline{on the  torus and on the sphere.  }}

 \vspace{2mm}

 {\bf 4.1. Quadratically integrable geodesic flows on the torus. }

Let  $L$ be a positive number.
 Denote by  $S_L$  the circle with a smooth parameter  $t\in R(\mbox{ mod } L)$.

 \begin{Df} A metric  $G$ on the torus  $T^2$ is called Liouville
 if for an appropriate positive  $L$ and for appropriate
 nonconstant positive functions  $f:S_1\to R$, $h:S_L\to R$ there exist
 a diffeomorphism  $\chi:T^2\to S_1\times S_L$ that takes the metric
  $G$ to the metric  $(f(x)+h(y))(dx^2+dy^2)$, where $x,y$ are
  parameters on  $S_1,S_L$, respectively.
  \end{Df}

 \begin{Df}   A metric  $G$ on the torus  $T^2$ is called pseudo-liouville
 if there exists a Liouville metric  $G_{liuv}$ on the torus  $T^2$ and
 a covering  $\rho:T^2\to T^2$ such that  $\rho^*(G_{liuv})=G$. \end{Df}

 {\bf Remark}  {\small Since the identity mapping $id:T^2\to T^2$  is a 1-sheet covering
  a Liouville metric is a pseudo-liouville metric. }

 \begin{Th}[\cite{BN, Mat1}]
 The geodesic flow of a metric  $G$ on the torus  $T^2$ is quadratically
 integrable iff the metric  $G$ is pseudo-liouville.
 \end{Th}

 \vspace{.2cm}

 {\bf 4.2. Quadratically integrable geodesic flows on the sphere.}
Consider the torus  $S_1\times S_L$. Let for smooth functions
$f:S_1\to    R$,
$h:S_L\to    R$
 the following conditions hold
 \begin{itemize} {\small
 \item[1] The functions  $f$,  $h$ are nonnegative.
 The  function $f$ equals
  zero only in the points  0 and
  $\frac{1}{2}$;

 The  function $h$ equals
  zero only in the points  0 and
  $\frac{L}{2}$;

 \item[2] $f''(0)= f''(\frac{1}{2})=
 h''(0)=h''(\frac{L}{2})\ne0$.

\item[3]  For any  $x$,   $y$ \ \ $f(x)=f(-x)$ and
 $h(y)=h(-y)$.
                                      }
  \end{itemize}

Consider  the involution  $\sigma:S_1\times S_L\to S_1\times S_L$,
 $\sigma(x,y)\eqdef(-x,-y)$. The fixed points of $\sigma$ are:
  $(0,0)$, $(0,\frac{L}{2})$,
 $(\frac{1}{2},0)$ and $(\frac{1}{2},\frac{L}{2})$.

Consider the factorspace
 $\tilde S\eqdef S_1\times S_L/_\sigma$.
 $\tilde S$ is homeomorphic to
 the sphere. Consider the following structure of a smooth  2-dimensional
 manifold on  $\tilde S$.

 Denote by  $\chi:S_1\times S_L\to \tilde S$
 the dual  to the involution
  $\tilde S$ mapping. In other words, the mapping $\chi$
  takes a point
                          $(x,y)\in S_1\times S_L$ to the point
 $((x,y),(-x,-y))\in \tilde S$.

 Consider smooth structure on  $\tilde S$ such that the mapping
  $\chi$  is a smooth branched covering
  with branch points
   $\chi(0,0)$
 $\chi(0,\frac{L}{2})$
 $\chi(\frac{1}{2},\frac{L}{2})$
 $\chi(\frac{1}{2},0)$ of the branch
  index 1.

 Obviously, the smooth structure exists.

Note, that the mapping  $\chi$
in the appropriate coordinates is
 the Weierstrass $\wp$-function with  half-periods  $(\fr{1}{2},0)$ and
 $(0,\fr{L}{2})$.

Consider the degenerate metric
 $((f(x)+h(y))(dx^2+dy^2)$     on the torus
  $S_1\times S_L$. The form
 $((f(x)+h(y))(dx^2+dy^2)$ is positive definite everywhere exept of at  the
 points
 $(0,0)$,
 $(0,\frac{L}{2})$,
 $(\frac{1}{2},\frac{L}{2})$,    and
 $(\frac{1}{2},0)$.
  Since the function $f(x)+h(y)$ is preserved   by the involution
   $\sigma$, we see that the metric  $((f(x)+h(y))(dx^2+dy^2)$
   induces a metric on the sphere
    $\tilde S$ without branch points.
     The conditions on the functions  $f$, $h$ allow to complement the
     metric in the branch points. Denote by
      $\tilde G_{L,f,h}$ the complemented metric.

 \begin{Df} A metric  $G$ on the sphere  $S^2$ is called a Kolokoltzov metric
 if for an appropriate positive number  $L>0$ and for appropriate
 functions  $f$, $h$  there exists a diffeomorphism
 $\tilde S\to S^2$ that takes the metric $\tilde G_{L,f,h}$ to the metric
  $G$.
 \end{Df}

 \begin{Th}[ \cite{Kol}] The geodesic flow of a metric
 $G$ on the sphere  $S^2$ is quadratically integrable
 iff the metric  $G$ is a Kolokoltzov metric.
                              \end{Th}

   {\bf Remark} {\small If the functions  $f$ and  $h$ are smooth functions,
   then the metric
 $\tilde  G_{L,f,h}$ is smooth on $\tilde S$ without the branch points.
 The metric  $\tilde G_{L,f,h}$
 is $C^k$-smooth in a branch point iff the following condition holds.
 In the branch point Taylor series of the function $f$ as the function of $x$
 coincides till  $(k+2)$-member
 with Taylor series of the  function $h$ as the function of $-x$ (see \cite{Kol2}).
 In other words, the metric
 $\tilde G_{L,f,h}$ is $C^k$-smooth in a branch point $(0,0)$
 if for any natural $m\le k+2$ \
 $\frac{d^m f}{d x^m}_{|_0}=(-1)^m\frac{d^m h}{d y^m}_{|_0}$.
 }

                \vspace{.2cm}

{\bf 4.3. Hyperbolic trajectories of the quadratically integrable
geodesic flows on the torus. }

Consider the isoenergy surface
$\{H=1\}=\{\frac{p_x^2+p_y^2}{f(x)+h(y)}=1\}.$

Following the paper  \cite{Sel} we describe the set of critical circles of
the geodesic flow. For simplicity suppose  $f,h$ are Morse functions.

Denote by  $K(f)$ ( $K(h)$ ) the set of critical points of the function
 $f$        (respetively, $h$).
 In  \cite{Sel}, E. Selivanova  proved that the set of critical points
 of the metric  $G$ is the union of the sets
  $${\cal
O}^+_f=\{(x,y,p_x,p_y):x\in K(f),p_x=0,p_y=\sqrt{f(x)+h(y)}\},$$
 $${\cal
O}^-_f=\{(x,y,p_x,p_y):x\in K(f),p_x=0,p_y=-\sqrt{f(x)+h(y)}\},$$
$${\cal
O}^+_h=\{(x,y,p_x,p_y): y\in K(h),p_x=\sqrt{f(x)+h(y)},p_y=0\},$$
 $${\cal
  O}^-_h=\{(x,y,p_x,p_y): y\in K(h),p_x=-\sqrt{f(x)+h(y)},p_y=0\}.$$
Besides, the circle
$$\{(x,y,p_x,p_y):x=x_0\in K(f),p_x=0,p_y=-\sqrt{f(x)+h(y)}\}$$ is a saddle
circle iff the point  $x_0$ is a nodegenerate critical  point of the
Morse index  $1$.

It is possible to prove that the saddle circles of quadratically
integrable geodesic
flow on the torus are hyperbolic
    trajectories.

     \vspace{.2cm}
 {\bf 4.4. Saddle circles of the quadratically integrable geodesic flow on the
 sphere.}

Suppose the functions  $f$,  $h$  satisfy the conditions
 1,2,3 from  section 4.2. Consider the torus
  $S_1\times S_L$ with the (degenerate) metric
$ds^2=(f(x)+h(y))(dx^2+dy^2).$  Consider the torus without
the points
$(0,0)$,
 $(0,\frac{L}{2})$,
 $(\frac{1}{2},\frac{L}{2})$,
 $(\frac{1}{2},0)$. Then   $ds^2$ is a metric.
 Using  the previous section,   we have that the saddle circles of the
 metric are the circles
 $$\{(x,y,p_x,p_y):x=x_0,p_x=0,p_y=\pm\sqrt{f(x)+h(y)}, \mbox{where  $x_0$
is a critical  point of the Morse
index 1}\}$$ and
$$\{(x,y,p_x,p_y):y=y_0,p_y=0,p_x=\pm\sqrt{f(x)+h(y)}, \mbox{where
 $y_0$
is a critical  point of the Morse index  1}\}.$$

The involution  $\sigma$ rerrange
pairs of circles. Since that, any
 pair factorize to a saddle circle. Such saddle circles are called
 {\it simple}. If the projection of the saddle circle  of the geodesic flow
 does not contain the point from   the set $\{$
$(0,0)$,
 $(0,\frac{L}{2})$,
 $(\frac{1}{2},\frac{L}{2})$,
 $(\frac{1}{2},0)\}$,         then  the circle is simple.

       There exist two  nonsimple saddle
       circles
(we denote them by $\hat\ga_1$,
$\hat\ga_2$)                   of the geodesic
flow of a Kolokoltzov metric.
Consider the segments
$$I_f^+=\{(x,y,p_x,p_y):x\in(0,\fr{1}{2}),y=0,p_y=0,p_x=\sqrt{f(x)+h(y)}\},$$
$$I_f^-=\{(x,y,p_x,p_y):x\in(0,\fr{1}{2}),y=\fr{1}{2},p_y=0,p_x=-\sqrt{f(x)+h(y)}\},$$
$$I_h^+=\{(x,y,p_x,p_y):x=0,y\in(0,\fr{L}{2}),p_y=0,p_y=\sqrt{f(x)+h(y)}\},$$
$$I_h^-=\{(x,y,p_x,p_y):x=0,y\in(0,\fr{L}{2}),p_y=0,p_y=-\sqrt{f(x)+h(y)}\}.$$
The segments factorize to the circle  ( denoted
 by $\hat\ga_1$),
Since the restriction of the integral to the transversal disk  to
an  intervals has a  singularity of Morse
index 1, we see that
$\hat\ga_1$ is a saddle circle.

Note, that the simple saddle circles of the geodesic flow
of a Kolokoltzov metric are hyperbolic trajectories.
There exist examples of the Kolokoltzov metric such that
the saddle circles
($\hat\ga_1$, $\hat\ga_2$) are not hyperbolic trajectories.

 \vspace{.2cm}
 {\sc \S5. Fundamental solution of the Jacobi equation for the hyperbolic geodesics
 of quadratically integrable geodesic flows.}

\vspace{.2cm}

 The formulas (3--5) allows to construct a fundamental solution of
 the Jacobi-Hill equation for saddle circles of the quadratically
 integrable geodesic flows.

{\bf 5.1. Torus. }
Without loss of generality it can be assumed that
a saddle circle as the set of points coincides with the set
$\{(x,y,p_x,p_y):y=0,p_y=0,p_x=\sqrt{f(x)+h(y)}\}$, and that
 $h(0)=0$.
Consider the Liouville fiber, which contains the saddle circle.
The Liouville fiber is the set
$\{(x,y,p_x,p_y):p_x=\sqrt{f(x)},p_y=\pm\sqrt{h(y)}\}$.
It is easy to see that the vectors
$Y_\pm=\p_x\pm\fr{\sqrt{f(x)}}{\sqrt{2}}\p_{p_x}$ are {\small \begin{itemize}
\item[{\bf 1)}] tangential to
the Liouville fiber
and \item[{\bf 2)}] lie in the planes  $<D_2 D_\phi>$.
\end{itemize}}
Therefore the vectors
$Y_\pm=\p_x\pm\fr{\sqrt{f(x)}}{\sqrt{2}}\p_{p_x}$  are direction vectors
of two families  of invariant 1-subspases.

Using Lemma~\ref{lemma}, we get the following equations for
 $\al$, $\beta$.

\begin{equation}\label{a}
\al(x)=\pm\sqrt{f(x)}
 \end{equation}
\begin{equation}\label{b}
\beta(x)=\fr{f'(x)}{2\sqrt{h}}-\sqrt{\fr{h''(0)}{2h}}
 \end{equation}

Combining  (\ref{a}, \ref{b}) with  the formula for  $u$,
we get a fundamental solution
of the Jacobi-Hill equation:
\begin{equation}
\label{u}
\matrix{u_+(x)=\sqrt{f(x)}\exp\left[\ \
\sqrt{\frac{h''(0)}{2}}\int\limits_{x_0}^x\fr{ds}{\sqrt{f(s)}}\right]\cr \
u_-(x)=\sqrt{f(x)}\exp\left[-\sqrt{\frac{h''(0)}{2}}\int\limits_{x_0}^x\fr{ds}{\sqrt{f(s)}}\right]\cr}
\end{equation}

\vspace{.2cm}
{\bf Remark } {\small  Since the function $u_+$  increase and since
the function $u_-$
 decrease, we see that there are no conjugate points
along  the hyperbolic geodesics
 of a quadratically integrable geodesic flow on the torus.

}

\vspace{.2cm}

 {\bf 5.2. Sphere.} For a simple geodesic
 the answer coinsides with the answer for
 the torus.

 Consider the nonsimple saddle circle  $\hat\ga_1$.
 The circle can be represented as four glued segments
  $I_h^+$, $I_f^+$,  $I_h^-$,  $I_f^-$. Using (\ref{u}), we see that
  the pair of functions
 \begin{equation}
\label{u1}
\matrix{u_+^f(x)=\sqrt{f(x)}\exp\left[\ \
\sqrt{\frac{h''(0)}{2}}\int\limits_{x_0}^x\fr{ds}{\sqrt{f(s)}}\right]\cr \
u_-^f(x)=\sqrt{f(x)}\exp\left[-\sqrt{\frac{h''(0)}{2}}
 \int\limits_{x_0}^x\fr{ds}{\sqrt{f(s)}}.\right]\cr}
\end{equation}
 is a  fundamental solution on the segments
  $I_f^+$ ¨ $I_f^-$.

   Arguing as above,
the pair of functions
\begin{equation}
\label{u2}
\matrix{u_+^h(y)=\sqrt{h(y)}\exp\left[\ \
\sqrt{\frac{f''(0)}{2}}\int\limits_{y_0}^y\fr{ds}{\sqrt{h(s)}}\right]\cr \
u_-^h(y)=\sqrt{h(y)}\exp\left[-\sqrt{\frac{f''
 (0)}{2}}\int\limits_{y_0}^y\fr{ds}{\sqrt{h(s)}}.\right]\cr}
\end{equation}
   is a fundamental solution on the segments  $I_h^+$, $I_h^-$.
   We  shall glue the fundamental solutions of the Jacobi-Hill eqution
   in the points
$(0,0)$,
 $(0,\frac{L}{2})$,
 $(\frac{1}{2},\frac{L}{2})$ ¨
 $(\frac{1}{2},0)$.
 Consider the point  $(0,0)$. We have to find constants
 $C_{11}$,
$C_{12}$,
$C_{21}$,
$C_{22}$ such that
$$\mbox{lim}_{x\to 0}C_{11}\left(\matrix{u_+^f\cr\dot
u_+^f\cr}\right)+C_{12}\left(\matrix{u_-^f\cr\dot u_-^f\cr}\right)=
\mbox{lim}_{y\to 0}\left(\matrix{u_+^h\cr\dot u_+^h\cr}\right),$$
$$\mbox{lim}_{x\to 0}C_{21}\left(\matrix{u_+^f\cr\dot
u_+^f\cr}\right)+C_{22}\left(\matrix{u_-^f\cr\dot u_-^f\cr}\right)=
\mbox{lim}_{y\to 0}\left(\matrix{u_-^h\cr\dot u_-^h\cr}\right).$$

Therefore,
\begin{equation}\label{c}
C_{11}=\fr{
\lim\limits_{y\to 0}
\sqrt{h(y)}\exp\left[\ \
\sqrt{\frac{f''(0)}{2}}\int\limits_{y_0}^y\fr{ds}{\sqrt{h(s)}}\right]
}
{
\lim\limits_{x\to 0}
\sqrt{f(x)}\exp\left[\ \
\sqrt{\frac{h''(0)}{2}}\int\limits_{x_0}^x\fr{ds}{\sqrt{f(s)}}\right]
}, \ \  C_{12}=0 \end{equation}
\begin{equation}\label{c1}
C_{21}=0, \ \ C_{22}=\fr{
\lim\limits_{x\to 0}
\sqrt{f(x)}\exp\left[\ \
\sqrt{\frac{h''(0)}{2}}\int\limits_{x_0}^x\fr{ds}{\sqrt{f(s)}}\right]
}
{
\lim\limits_{y\to 0}
\sqrt{h(y)}\exp\left[\ \
\sqrt{\frac{f''(0)}{2}}\int\limits_{y_0}^y\fr{ds}{\sqrt{h(s)}}\right]
}.\end{equation}

Arguing as above,  we can
glue the fundamential solutions  in the points
$(0,0)=A_1$,
$(\fr{1}{2},\fr{L}{2})=A_2$,
 $(\fr{1}{2},0)=B_1$,
$(0,\fr{L}{2})=B_2$. We obtain a fundamental solution
$ (x_+, x_-)$ of the Jacobi-Hill equation along the geodesic line
 $\pi(\ga_1)$.




Consider  a point
$x_1\in I_f^+$.  Using the Sturm-Liouville theorem, we see that there exists a
point
$(\fr{1}{2}-x_2)\in I_f^-$ that is conjugate to the point $x_1$.
Using  (\ref{c}-\ref{c1}), we get the following equation for $x_2$:
\begin{equation}\label{con}
\matrix{{\lim\limits_{y\to 0}\sqrt{h(y)h(\fr{L}{2}-y})\exp\left[\ \
\sqrt{\frac{f''(0)}{2}}\int\limits_{y}^{\fr{L}{2}-y}\fr{ds}{\sqrt{h(s)}}\right]
}
=\cr
=\fr{
\lim\limits_{x\to 0}
\sqrt{f(x)}\exp\left[\ \
\sqrt{\frac{h''(0)}{2}}\int\limits_{x_1}^{\fr{1}{2}-x}
\fr{ds}{\sqrt{f(s)}}\right]
}{
\lim\limits_{x\to 0}
\sqrt{f(x)}\exp\left[\ \
\sqrt{\frac{h''(0)}{2}}\int\limits_{x_2}^x\fr{ds}{\sqrt{f(s)}}\right]
}}                         \end{equation}

\end{document}